\documentclass[11pt,a4paper]{article}

\usepackage{inputenc}
\usepackage{amsmath}
\usepackage{bm}
\usepackage{bbold}
\usepackage{amsthm}
\usepackage{enumerate}

\usepackage{hyperref}

\title{Tropical Optimization Problems\thanks{Advances in Economics and Optimization: Collected Scientific Studies Dedicated to the Memory of L. V. Kantorovich (L. A. Petrosyan, J. V. Romanovsky, D. W. K. Yeung, eds), pp. 195-214, Nova Science Publishers, New York, 2014.}}

\author{Nikolai Krivulin\thanks{Faculty of Mathematics and Mechanics, Saint Petersburg State University, 28 Universitetsky Ave., Saint Petersburg, 198504, Russia, 
nkk@math.spbu.ru.}\thanks{The work was supported in part by the Russian Foundation for Humanities under Grant \#13-02-00338.}
}

\date{}

\newtheorem{theorem}{Theorem}
\newtheorem{lemma}[theorem]{Lemma}
\newtheorem{corollary}[theorem]{Corollary}

\begin{document}

\maketitle

\begin{abstract}
We consider optimization problems that are formulated and solved in the framework of tropical mathematics. The problems consist in minimizing or maximizing functionals defined on vectors of finite-dimen\-sional semimodules over idempotent semifields, and may involve constraints in the form of linear equations and inequalities. The objective function can be either a linear function or nonlinear function calculated by means of multiplicative conjugate transposition of vectors. We start with an overview of known tropical optimization problems and solution methods. Then, we formulate certain new problems and present direct solutions to the problems in a closed compact vector form suitable for further analysis and applications. For many problems, the results obtained are complete solutions.
\\

\textbf{Key-Words:} idempotent semifield, tropical optimization, nonlinear objective function, linear constraints, Chebyshev approximation.
\\

\textbf{MSC (2010):} 65K10, 15A80, 65K05, 90C48, 41A50
\end{abstract}

\section{Introduction}

Tropical (idempotent) mathematics presents an area of applied mathematics that focuses on the theory and applications of semirings with idempotent addition. Among the first works, which laid the foundation for the area, were the pioneering papers by Pandit \cite{Pandit1961Anew}, Cuninghame-Green \cite{Cuninghamegreen1962Describing}, Giffler \cite{Giffler1963Scheduling} and Hoffman \cite{Hoffman1963Onabstract}, dating back to the early 1960s. 

The classical works published at the same time or somewhat earlier by Kleene \cite{Kleene1956Representation}, Bellman \cite{Bellman1956Onaquasilinear}, and Bellman and Karush \cite{Bellman1961Onanewfunctional,Bellman1962Mathematical} also served as an essential prerequisite. Specifically, a significant effect was exerted on the area by the Bellman--Karush maximum transform. This transform, which was initially introduced in \cite{Bellman1961Onanewfunctional} in a multiplicative form, has taken an additive form in \cite{Bellman1962Mathematical} and in many later publications. The additive form of the transform is in direct connection with the theory of support functions in convex optimization \cite{Rockafellar1973Convex}. Note that this additive form was also proposed by Romanovsky \cite{Romanovsky1962Someremarks}.

A notable contribution to the development of the area was made by the Leningrad (St.~Petersburg) researchers N.~N.~Vorob{'}ev, J.~V.~Romanovsky and A.~A.~Korbut. In the paper \cite{Vorobjev1963Theextremal} published in 1963, Vorob{'}ev provided a groundwork for the theory of linear operators on extremal spaces, which were defined as vector spaces over idempotent semifields with the operations of maximum and minimum in the role of addition.                                                                       

Then came the works \cite{Romanovskii1964Asymptotic,Romanovskii1965Asymptotic,Romanovskii1967Optimization,Romanovskiy1972Deterministic}, in which Romanovsky studied asymptotic properties of solutions in dynamic programming problems and obtained an idempotent analogue for the Gelfand spectral radius theorem for bounded linear operators in Banach spaces. These results have proved to be conceptually close to the theory of the extremal spaces, which was further developed by Vorob{'}ev \cite{Vorobjev1967Extremal,Vorobjev1970Extremal} and Korbut \cite{Korbut1965Extremal,Korbut1972Extremal}.

An extensive advancement in the theory and applications of tropical mathematics was achieved by V.~P.~Maslov and his colleagues, who have introduced idempotent functional analysis as a new important field within the area. The results obtained by the Maslov scientific school and published in the monographs by Maslov and Samborski{\u\i} \cite{Maslov1992Idempotent}, Maslov and Kolokoltsov \cite{Maslov1994Idempotent}, Kolokoltsov and Maslov \cite{Kolokoltsov1997Idempotent}, as well as in many research papers, have provided a strong analytical framework for the contemporary tropical mathematics. 

In the last few decades, tropical mathematics has received considerable attention as both a powerful mathematics apparatus and a flexible application tool. Advances in the theory and applications have resulted in a great many publications, including the monographs by Cuninghame-Green \cite{Cuninghamegreen1979Minimax}, Carr\'{e} \cite{Carre1979Graphs}, U.~Zimmermann \cite{Zimmermann1981Linear}, Baccelli et al. \cite{Baccelli1993Synchronization}, Golan \cite{Golan2003Semirings}, Heidergott, Olsder and van der Woude \cite{Heidergott2006Maxplus}, Gondran and Minoux \cite{Gondran2008Graphs}, Krivulin \cite{Krivulin2009Methods}, Butkovi\v{c} \cite{Butkovic2010Maxlinear}, McEneaney \cite{Mceneaney2010Maxplus}, and in a large amount of contributed papers.

Optimization problems, which can be formulated and solved in terms of tropical mathematics, form an important research domain within the area. These problems arise in many real-world applications in job scheduling, location analysis, transportation networks, decision making, discrete event dynamic systems, and in other fields.

Optimization problems have been examined even in the early publications \cite{Cuninghamegreen1962Describing,Hoffman1963Onabstract}. Since that time a number of studies that investigate optimization problems in the context of tropical mathematics has been published, which includes the monographs \cite{Cuninghamegreen1979Minimax,Zimmermann1981Linear,Krivulin2009Methods,Butkovic2010Maxlinear}.

Tropical optimization problems are commonly formulated to minimize or maximize functionals defined on vectors of a finite-dimensional semimodule over an idempotent semifield, and may include additional constraints in the form of equations and inequalities. Some problems have both the objective function and constraints that are linear in the sense of tropical mathematics. These problems, which can be considered as direct analogues to linear programming problems, were examined in \cite{Hoffman1963Onabstract,Superville1978Various,Zimmermann1981Linear}.

In another class of problems \cite{Cuninghamegreen1976Projections,Cuninghamegreen1979Minimax,Zimmermann1981Linear}, the constraints are linear, but the objective functions are not. In these problems, the objective function can usually be written by means of a multiplicative conjugate transposition operator, which does not preserve linearity.

To solve the problems considered, a variety of methods were proposed. Some problems in a rather general setting can be completely solved in an explicit form \cite{Hoffman1963Onabstract,Cuninghamegreen1976Projections,Superville1978Various,Zimmermann1981Linear}. For other problems, only algorithmic solutions are known in the form of iterative numerical procedures that produce a particular solution or indicate that no solution exists \cite{Butkovic1984Onproperties,Zimmermann1984Some}.

In this chapter, we offer a brief overview of known tropical optimization problems and solution methods. We formulate certain new problems and present direct solutions to the problems in a closed compact vector form suitable for further analysis and applications. For many problems, the results obtained are complete solutions.

The chapter is organized as follows. In Section~\ref{S-PDNR}, we briefly outline the definitions, notation, and key facts of tropical mathematics to be used in the subsequent sections. Section~\ref{S-OKPSM} is devoted to known optimization problems. We present new tropical optimization problems and describe obtained solutions in Section~\ref{S-SNPRR}. The chapter concludes with an extensive list of references on the topic considered.

\section{Preliminary Definitions, Notation and Results}
\label{S-PDNR}

We start with a brief overview of the basic elements of tropical (idempotent) mathematics to provide a formal unifying framework for formulating tropical optimization problems and representing the solutions to the problems in the subsequent sections. The overview mainly follows the notation system and results developed by Krivulin \cite{Krivulin2005Evaluation,Krivulin2006Eigenvalues,Krivulin2006Solution,Krivulin2009Methods}, 
which underlie the solutions presented. 

The literature on the topic offers a wide range of publications written in varied styles and available on both the introductory and advanced levels. Specifically, the recent monographs: Kolokoltsov and Maslov \cite{Kolokoltsov1997Idempotent}, Golan \cite{Golan2003Semirings}, Heidergott, Olsder and van der Woude \cite{Heidergott2006Maxplus}, Gondran and Minoux \cite{Gondran2008Graphs,Butkovic2010Maxlinear}, Butkovi\v{c} \cite{Butkovic2010Maxlinear}; and overview papers: Kolokoltsov \cite{Kolokoltsov2001Idempotent}, Litvinov \cite{Litvinov2007Themaslov}, Akian, Bapat and Gaubert \cite{Akian2007Maxplus}, provide a broad source for further reading.

\subsection{Idempotent Semifield}

An idempotent semifield is an algebraic system $(\mathbb{X},\oplus,\otimes,\mathbb{0},\mathbb{1})$, where $\mathbb{X}$ is a nonempty set closed under addition $\oplus$ and multiplication $\otimes$, and equipped with zero $\mathbb{0}$ and one $\mathbb{1}$, such that $(\mathbb{X},\oplus,\mathbb{0})$ is a commutative idempotent monoid, $(\mathbb{X},\otimes,\mathbb{1})$ is an abelian group, multiplication distributes over addition, and $\mathbb{0}$ is absorbing for multiplication.

In addition, we assume the semifield to admit a linear order, which is consistent with the partial order induced by idempotent addition to define $x\leq y$ if and only if $x\oplus y=y$. Moreover, the semifield is assumed to be algebraically complete (radicable), which implies that the equation $x^{p}=a$ has a root for all $a\in\mathbb{X}$ and every natural $p$.

Examples of the linearly ordered and algebraically complete idempotent semifield include 
$\mathbb{R}_{\max,+}=(\mathbb{R}\cup\{-\infty\},\max,+,-\infty,0)$, 
$\mathbb{R}_{\min,+}=(\mathbb{R}\cup\{+\infty\},\min,+,+\infty,0)$, 
$\mathbb{R}_{\max,\times}=(\mathbb{R}_{+}\cup\{0\},\max,\times,0,1)$ and 
$\mathbb{R}_{\min,\times}=(\mathbb{R}_{+}\cup\{+\infty\},\min,\times,+\infty,1)$, 
where $\mathbb{R}$ is the set of real numbers and $\mathbb{R}_{+}=\{x\in\mathbb{R}|x>0\}$.

Below, we drop the multiplication sign $\otimes$ for the brevity sake. The relation symbols and the optimization 
problems are thought of in terms of the linear order on the semifield.

\subsection{Algebra of Matrices and Vectors}

The set of matrices with $m$ rows and $n$ columns over $\mathbb{X}$ is denoted by $\mathbb{X}^{m\times n}$. A matrix with all entries equal to $\mathbb{0}$ is the zero matrix. A matrix is column-regular (row-regular) if it has no zero columns (rows). A matrix that has no zero rows nor zero columns is called regular.

Addition and multiplication of conforming matrices and scalar multiplication are performed component-wise in the usual way, except that the ordinary scalar addition is replaced by $\oplus$ and the multiplication by $\otimes$.

The multiplicative conjugate transpose of any nonzero matrix $\bm{A}=(a_{ij})\in\mathbb{X}^{m\times n}$ is the matrix $\bm{A}^{-}=(a_{ij}^{-})\in\mathbb{X}^{n\times m}$, where $a_{ij}^{-}=a_{ji}^{-1}$ if $a_{ji}\ne\mathbb{0}$, and $a_{ij}^{-}=\mathbb{0}$ otherwise.

A matrix with only one column (row) is a column (row) vector. The set of column vectors of size $n$ forms an idempotent semimodule over $\mathbb{X}$ and is denoted $\mathbb{X}^{n}$. A vector with all elements equal to $\mathbb{0}$ is the zero vector. A vector that has no zero elements is regular.

The multiplicative conjugate transpose of any nonzero vector $\bm{x}=(x_{i})\in\mathbb{X}^{n}$ is the row vector $\bm{x}^{-}=(x_{i}^{-})$, where $x_{i}^{-}=x_{i}^{-1}$ if $x_{i}\ne\mathbb{0}$, and $x_{i}^{-}=\mathbb{0}$ otherwise.

For any vector $\bm{x}\in\mathbb{X}^{n}$ and matrix $\bm{A}\in\mathbb{X}^{m\times n}$, idempotent analogues of the vector and matrix norms are given by
$$
\|\bm{x}\|
=
\bigoplus_{i=1}^{n}x_{i},
\qquad
\|\bm{A}\|
=
\bigoplus_{i=1}^{m}\bigoplus_{j=1}^{n}a_{ij}.
$$

With the vector $\mathbb{1}=(\mathbb{1},\ldots,\mathbb{1})^{T}$, which has all components equal to $\mathbb{1}$, we write
$$
\|\bm{x}\|
=
\mathbb{1}^{T}\bm{x},
\qquad
\|\bm{A}\|
=
\mathbb{1}^{T}\bm{A}\,\mathbb{1}.  
$$

Consider square matrices of order $n$ in $\mathbb{X}^{n\times n}$. A matrix with $\mathbb{1}$'s on the diagonal and $\mathbb{0}$'s elsewhere is the identity matrix denoted $\bm{I}$.

The trace of any matrix $\bm{A}=(a_{ij})$ is given by $\mathop\mathrm{tr}\bm{A}=a_{11}\oplus\cdots\oplus a_{nn}$.

To represent solutions to optimization problems below, we use a function that takes any matrix $\bm{A}$ to the scalar
$$
\mathop\mathrm{Tr}(\bm{A})
=
\mathop\mathrm{tr}\bm{A}\oplus\cdots\oplus\mathop\mathrm{tr}\bm{A}^{n}.
$$

If the condition $\mathop\mathrm{Tr}(\bm{A})\leq\mathbb{1}$ holds, we apply the asterate operator, 
which produces the matrix (sometimes called the Kleene star)
$$
\bm{A}^{\ast}
=
\bm{I}\oplus\bm{A}\oplus\cdots\oplus\bm{A}^{n-1}.
$$

If $\mathop\mathrm{Tr}(\bm{A})=\mathbb{1}$, we form the matrix $\bm{A}^{+}$, which holds those columns 
of the matrix $\bm{A}^{\times}=\bm{A}\bm{A}^{\ast}$ that have $\mathbb{1}$ on the diagonal.

\subsection{Linear Operators and Their Spectral Properties}

Every matrix $\bm{A}\in\mathbb{X}^{n\times n}$ defines a linear operator on $\mathbb{X}^{n}$, which is endowed with certain spectral properties. A scalar $\lambda\in\mathbb{X}$ and a nonzero vector $\bm{x}\in\mathbb{X}^{n}$ are corresponding eigenvalue and eigenvector of $\bm{A}$ if they satisfy the equation
$$
\bm{A}\bm{x}
=
\lambda\bm{x}.
$$

The matrix $\bm{A}$ always possesses an eigenvalue that is calculated by the formula
$$
\lambda
=
\bigoplus_{m=1}^{n}\mathop\mathrm{tr}\nolimits^{1/m}(\bm{A}^{m}).
$$

Any other eigenvalue of $\bm{A}$, if it exists, is not greater (in the sense of the order on $\mathbb{X}$) than $\lambda$, which is called the spectral radius of $\bm{A}$.

All eigenvectors of the matrix $\bm{A}$, which correspond to the spectral radius $\lambda$, are given by
$$
\bm{x}
=
(\lambda^{-1}\bm{A})^{+}\bm{u},
$$
where $\bm{u}$ is any regular vector of appropriate size.

\section{An Overview of Known Problems and Solution Methods}
\label{S-OKPSM}

Multidimensional optimization problems form an important research domain in tropical mathematics, dating back to the early works by Cuninghame-Green \cite{Cuninghamegreen1962Describing} and Hoffman \cite{Hoffman1963Onabstract}. These problems appear in various application areas, including job scheduling \cite{Cuninghamegreen1962Describing,Cuninghamegreen1976Projections,Cuninghamegreen1979Minimax,Zimmermann1981Linear,Zimmermann1984Some,Zimmermann2006Interval,Butkovic2009Introduction,Butkovic2009Onsome,Tam2010Optimizing,Aminu2012Nonlinear}, location analysis \cite{Cuninghamegreen1991Minimax,Krivulin2011Anextremal,Krivulin2012Anew}, transportation networks \cite{Zimmermann1981Linear,Zimmermann2006Interval}, decision making \cite{Elsner2004Maxalgebra,Elsner2010Maxalgebra,Gursoy2013Theanalytic}, and discrete event dynamic systems \cite{Gaubert1995Resource,Deschutter1996Maxalgebraic,Deschutter2001Model,Krivulin2005Evaluation}.

In this section, we offer a brief overview of known tropical optimization problems and discuss existing solution methods. We consider problems of minimizing or maximizing both linear and nonlinear functionals defined on vectors in finite-dimensional semimodules over idempotent semifields, subject to constraints expressed as linear equations and inequalities. The overview covers the problems with nonlinear objective functions that are, or can be, written by means of multiplicative conjugate transposition of vectors. Our analysis of the literature shows that many problems, which are relevant to tropical optimization, have objective functions that admit this form.   

Note that the problems with nonlinear objective functions are usually considered in a different way; these problems are formulated and solved either in the usual setting within the framework of ordinary mathematics \cite{Zimmermann1984Onmaxseparable,Zimmermann1984Some,Zimmermann1992Optimization,Zimmermann2003Disjunctive,Zimmermann2006Interval}, or in terms of two dual semifields $\mathbb{R}_{\max,+}$ and $\mathbb{R}_{\min,+}$ at once \cite{Cuninghamegreen1976Projections,Cuninghamegreen1979Minimax,Butkovic2009Onsome,Tam2010Optimizing}.

In the description of the problems below, we use the capital letters $\bm{A}$, $\bm{B}$ and $\bm{C}$ for known matrices and the lower case letters $\bm{b}$, $\bm{d}$, $\bm{p}$ and $\bm{q}$ for known vectors. The symbol $\bm{x}$ represents the unknown vector. The matrix and vector operations are thought of as defined in terms of an idempotent semifield. The minus sign in the exponent indicates multiplicative conjugate transposition of vectors and matrices.

Specifically, in terms of the real semifield $\mathbb{R}_{\max,+}$, vector addition, matrix-vector multiplication and scalar multiplication follow the standard rules, where the operations $\max$ and $+$ play the respective roles of addition and multiplication. Multiplicative conjugate transposition of a vector implies replacing each element of the vector by its inverse (which coincides with its opposite in ordinary arithmetic), and transposition.

\subsection{Problems with Linear Objective Functions}

One of the first optimization problems defined and solved in terms of tropical mathematics is a formal analogue of linear programming problems, which takes the form 
\begin{equation*}
\begin{aligned}
&
\text{minimize}
&&
\bm{p}^{T}\bm{x},
\\
&
\text{subject to}
&&
\bm{A}\bm{x}
\geq
\bm{d}.
\end{aligned}
\end{equation*}

Exact general solutions to the problem are obtained under various assumptions concerning the idempotent semiring, in which context the problem is examined. Specifically, Hoffman \cite{Hoffman1963Onabstract} considered a rather general idempotent semiring and proposed a solution based on abstract extension of the duality principle in linear programming. Another general solution was established by U.~Zimmermann \cite{Zimmermann1981Linear} by means of a residual-based solution technique. The common approach in \cite{Hoffman1963Onabstract} was further developed by Superville \cite{Superville1978Various} to examine the problem in the context of the semifield $\mathbb{R}_{\max,+}$. The problem was also solved by Gavalec and K.~Zimmermann \cite{Gavalec2012Duality} for the semifields $\mathbb{R}_{\max,+}$ and $\mathbb{R}_{\max,\times}$ in the framework of max-separable functions.

The results of the theory of max-separable functions were applied by K.~Zimmermann \cite{Zimmermann1984Onmaxseparable,Zimmermann1984Some,Zimmermann1992Optimization,Zimmermann2003Disjunctive,Zimmermann2006Interval} to solve the problem with more constraints 
\begin{equation*}
\begin{aligned}
&
\text{minimize}
&&
\bm{p}^{T}\bm{x},
\\
&
\text{subject to}
&&
\bm{A}\bm{x}
\leq
\bm{d},
\quad
\bm{C}\bm{x}
\geq
\bm{b},
\\
&&&
\bm{g}
\leq
\bm{x}
\leq
\bm{h}.
\end{aligned}
\end{equation*}

Under general conditions, an exact solution was obtained, which was, however, given in ordinary terms rather than in terms of tropical mathematics.

Butkovi\v{c} \cite{Butkovic1984Onproperties}, Butkovi\v{c} and Aminu \cite{Butkovic2009Introduction,Aminu2012Nonlinear} studied a problem that has constraints in the form of two-sided equation (with the unknown vector on both left and right sides), and can be written as
\begin{equation*}
\begin{aligned}
&
\text{minimize}
&&
\bm{p}^{T}\bm{x},
\\
&
\text{subject to}
&&
\bm{A}\bm{x}\oplus\bm{b}
=
\bm{C}\bm{x}\oplus\bm{d}.
\end{aligned}
\end{equation*}

Specifically, a pseudo-polynomial algorithm was proposed in \cite{Butkovic2009Introduction}, which produces a solution if solutions exist or indicates that the problem has no solutions. The algorithm applies an alternating method proposed in \cite{Cuninghamegreen2003Theequation} to replace the equality constraints by two opposite inequalities, which are alternately solved to achieve more and more accurate estimates for the solution. An heuristic approach that combines a search scheme providing approximate solutions with iterative procedures of solving low-dimensional problems with one or two unknown variables was developed in \cite{Aminu2012Nonlinear}.

\subsection{Problems with Nonlinear Objective Functions}

The tropical optimization problems with nonlinear objective functions, in which nonlinearity results from the application of multiplicative conjugate transposition, form a rich class of various problems that arise in many applications. We divide the problems into a few groups according to the particular form of the objective function and the principal interpretation of the problems.

\subsubsection{Chebyshev Approximation}

Among the first problems with nonlinear objective functions was the problem, which is formulated in the form
\begin{equation*}
\begin{aligned}
&
\text{minimize}
&&
(\bm{A}\bm{x})^{-}\bm{p},
\\
&
\text{subject to}
&&
\bm{A}\bm{x}
\leq
\bm{p}.
\end{aligned}
\end{equation*}

The problem was considered by Cuninghame-Green \cite{Cuninghamegreen1976Projections} in the context of approximation in the semifield $\mathbb{R}_{\max,+}$ with the Chebyshev metric. The problem required finding a vector $\bm{x}$ that gives the best underestimating approximation to a vector $\bm{p}$ by means of vectors $\bm{A}\bm{x}$. A general solution in an explicit form was obtained based on the theory of linear operators on vectors over the semifield $\mathbb{R}_{\max,+}$. A similar solution was provided by U.~Zimmermann \cite{Zimmermann1981Linear}.

An unconstrained approximation problem formulated in terms of a general idempotent semifield was examined by Krivulin \cite{Krivulin2009Methods} in the form
\begin{equation}
\begin{aligned}
&
\text{minimize}
&&
(\bm{A}\bm{x})^{-}\bm{p}\oplus\bm{p}^{-}\bm{A}\bm{x}.
\end{aligned}
\label{P-AxppAx}
\end{equation}

To solve the problem, an approach was proposed, which involves the derivation of a sharp lower bound for the objective function. The form of this bound is exploited to construct vectors that produce the bound. The results obtained were applied to establish existence conditions and to describe all solutions in a closed vector form for the equation $\bm{A}\bm{x}=\bm{p}$. In addition, the following problems of both underestimating and overestimating approximation were solved as consequences of the results: 
\begin{equation*}
\begin{aligned}
&
\text{minimize}
&&
(\bm{A}\bm{x})^{-}\bm{p},
\\
&
\text{subject to}
&&
\bm{A}\bm{x}
\leq
\bm{p};
\end{aligned}
\qquad\qquad\qquad
\begin{aligned}
&
\text{minimize}
&&
\bm{p}^{-}\bm{A}\bm{x},
\\
&
\text{subject to}
&&
\bm{A}\bm{x}
\geq
\bm{p}.
\end{aligned}
\end{equation*}

Certain of the problems, which were originally formulated and solved in a different setting, can be represented in terms of tropical mathematics as problems with nonlinear objective functions and linear constraints and thus are also included in the overview. In particular, a constrained problem of the Chebyshev approximation in the semifield $\mathbb{R}_{\max,+}$ was examined in \cite{Zimmermann1984Some}, where the problem was solved with a polynomial-time threshold-type algorithm. In terms of tropical mathematics, this problem can be written in the form
\begin{equation}
\begin{aligned}
&
\text{minimize}
&&
(\bm{A}\bm{x})^{-}\bm{p}\oplus\bm{p}^{-}\bm{A}\bm{x},
\\
&
\text{subject to}
&&
\bm{g}
\leq
\bm{x}
\leq
\bm{h}.
\end{aligned}
\label{P-AxppAxgxh}
\end{equation}

\subsubsection{Problems with Span Seminorm}

The problems examined by Butkovi\v{c} and Tam \cite{Butkovic2009Onsome,Tam2010Optimizing} in the context of the semifield $\mathbb{R}_{\max,+}$ had the objective function defined as the span seminorm, which calculates the maximum deviation between elements of a vector. Solutions to both minimization and maximization problems were obtained using a combined technique based on the representation in terms of both the semifield $\mathbb{R}_{\max,+}$ and its dual $\mathbb{R}_{\min,+}$. Note that these problems can be described only in terms of the semifield $\mathbb{R}_{\max,+}$ as follows
\begin{equation}
\begin{aligned}
&
\text{minimize}
&&
\mathbb{1}^{T}\bm{A}\bm{x}(\bm{A}\bm{x})^{-}\mathbb{1},
\end{aligned}
\qquad\qquad
\begin{aligned}
&
\text{maximize}
&&
\mathbb{1}^{T}\bm{A}\bm{x}(\bm{A}\bm{x})^{-}\mathbb{1},
\end{aligned}
\label{P-minmax1AxAx1}
\end{equation}
where $\mathbb{1}=(\mathbb{1},\ldots,\mathbb{1})^{T}$ is a vector of ones in the sense of $\mathbb{R}_{\max,+}$.

\subsubsection{A Problem of ``Linear-fractional'' Programming}

A minimization problem that has constraints in the form of a two-sided inequality, and that is written in terms of the semifield $\mathbb{R}_{\max,+}$ as
\begin{equation*}
\begin{aligned}
&
\text{minimize}
&&
\bm{p}^{T}\bm{x}(\bm{q}^{T}\bm{x})^{-1},
\\
&
\text{subject to}
&&
\bm{A}\bm{x}\oplus\bm{b}
\leq
\bm{C}\bm{x}\oplus\bm{d},
\end{aligned}
\end{equation*}
was investigated by Gaubert, Katz and Sergeev \cite{Gaubert2012Tropical}. The problem, called a tropical linear-fractional programming problem, was solved by means of an iterative computational scheme that exploits the relationship established by Akian, Gaubert and Guterman \cite{Akian2012Tropical} between the solutions to two-sided vector equations in the sense of $\mathbb{R}_{\max,+}$ and mean payoff games.

\subsubsection{Extremal Property of the Spectral Radius}

A problem examined by Cuninghame-Green \cite{Cuninghamegreen1962Describing} to minimize a functional defined on vectors over the semifield $\mathbb{R}_{\max,+}$ was apparently the first optimization problem, which arose in the context of tropical mathematics. With the use of multiplicative conjugate transposition, the problem can be represented in the form 
\begin{equation}
\begin{aligned}
&
\text{minimize}
&&
\bm{x}^{-}\bm{A}\bm{x}.
\end{aligned}
\label{P-minxAx}
\end{equation}

As one of the main result in \cite{Cuninghamegreen1962Describing}, which gave impetus to further research on the spectral theory of linear operators in tropical mathematics, it was shown that the minimum in the problem is equal to the spectral radius $\lambda$ of the matrix $\bm{A}$ and is attained at any eigenvector that corresponds to the radius. Moreover, explicit expressions were obtained to calculate the spectral radius and the eigenvector in terms of standard arithmetic operations. Similar results were provided by Engel and Schneider \cite{Engel1975Diagonal} and by Superville \cite{Superville1978Various}.

Furthermore, the above solution was extended and represented in more general terms of tropical mathematics by Cuninghame-Green \cite{Cuninghamegreen1979Minimax}. Analogues results in a compact vector form using multiplicative conjugate transposition were proposed by Krivulin \cite{Krivulin2005Evaluation,Krivulin2006Eigenvalues,Krivulin2009Methods}.

To find all vectors that yield the minimum in the problem, a computational approach was proposed in \cite{Cuninghamegreen1979Minimax}. It consists in solving a linear programming problem and thus does not guarantee a direct solution in an explicit form.

Finally, Elsner and van den Driessche \cite{Elsner2004Maxalgebra,Elsner2010Maxalgebra} indicated that the solution set includes, in addition to the eigenvectors, all vectors $\bm{x}$ that satisfy the inequality $\bm{A}\bm{x}\leq\lambda\bm{x}$. An iterative computational procedure was offered to obtain vectors from this solution set.

\section{Some New Problems and Recent Results}
\label{S-SNPRR}

In this part, we present several new problems 
\cite{Krivulin2014Aconstrained,Krivulin2013Amaximization,Krivulin2014Complete,Krivulin2013Direct,Krivulin2013Explicit,Krivulin2013Extremal}, which are formulated and solved in terms of a general algebraically complete linearly ordered idempotent semifield. If a complete solution to a problem is obtained, the solution is given below by an expression that represents a general vector describing all vectors in the solution set, and that is referred to as the general solution.

\subsection{Chebyshev-like Approximation Problems}

We start with constrained optimization problems with nonlinear objective functions that have the form similar to that in the Chebyshev approximation problems \eqref{P-AxppAx} and \eqref{P-AxppAxgxh}. One of the application areas for these problems is the single facility location problems in multidimensional spaces with Chebyshev distance under various constraints in the form of linear equalities and inequalities (see, e.~g., \cite{Krivulin2012Anew,Krivulin2013Direct,Krivulin2014Complete}). 

To solve the next two problems with simple boundary constraints, we apply a general approach which was proposed and developed in \cite{Krivulin2005Evaluation,Krivulin2006Eigenvalues,Krivulin2009Methods}. The approach involves the derivation of a sharp lower bound on the objective function and the use of this bound to obtain solution vectors.

First, consider the following problem: given vectors $\bm{p},\bm{q},\bm{g},\bm{h}\in\mathbb{X}^{n}$, find regular vectors $\bm{x}\in\mathbb{X}^{n}$ such that
\begin{equation}
\begin{aligned}
&
\text{minimize}
&&
\bm{q}^{-}\bm{x}\oplus\bm{x}^{-}\bm{p},
\\
&
\text{subject to}
&&
\bm{g}
\leq
\bm{x}
\leq
\bm{h}.
\end{aligned}
\label{P-qxxpgxh}
\end{equation}

\begin{theorem}[\cite{Krivulin2013Direct}]
\label{T-qxxpgxh}
Let $\bm{x}$ be the general regular solution to problem \eqref{P-qxxpgxh}, where $\bm{p}$ and $\bm{q}$ are regular vectors, $\bm{g}$ and $\bm{h}$ are vectors such that $\bm{g}\leq\bm{h}$, and let $\Delta=\sqrt{\bm{q}^{-}\bm{p}}$.

Then, the minimum in \eqref{P-qxxpgxh} is equal to
\begin{equation*}
\mu
=
\Delta\oplus\bm{q}^{-}\bm{g}\oplus\bm{h}^{-}\bm{p},
\end{equation*}
and the general solution is given by
\begin{equation*}
\mu^{-1}\bm{p}\oplus\bm{g}
\leq
\bm{x}
\leq
(\mu^{-1}\bm{q}^{-}\oplus\bm{h}^{-})^{-}.
\end{equation*}
\end{theorem}

Suppose that, given a matrix $\bm{A}\in\mathbb{X}^{m\times n}$ and vectors $\bm{p},\bm{q}\in\mathbb{X}^{m}$, $\bm{g}\in\mathbb{X}^{n}$, the problem is to find regular vectors $\bm{x}\in\mathbb{X}^{n}$ that
\begin{equation}
\begin{aligned}
&
\text{minimize}
&&
\bm{q}^{-}\bm{A}\bm{x}\oplus(\bm{A}\bm{x})^{-}\bm{p},
\\
&
\text{subject to}
&&
\bm{x}
\geq
\bm{g}.
\end{aligned}
\label{P-qAxAxpxg}
\end{equation}

\begin{theorem}[\cite{Krivulin2013Direct}]
\label{T-qAxAxpxg}
Let $\bm{x}$ be a regular solution to problem \eqref{P-qAxAxpxg}, where $\bm{A}$ is a regular matrix, $\bm{p}$ and $\bm{q}$ are regular vectors, and let $\Delta=\sqrt{(\bm{A}(\bm{q}^{-}\bm{A})^{-})^{-}\bm{p}}$.

Then, the minimum in \eqref{P-qAxAxpxg} is equal to
$$
\mu
=
\Delta\oplus\bm{q}^{-}\bm{A}\bm{g},
$$
and attained at
\begin{equation*}
\bm{x}
=
\mu(\bm{q}^{-}\bm{A})^{-}.
\end{equation*}
\end{theorem}

We now consider another problem with conditions that include a linear inequality with a matrix. Given vectors $\bm{p},\bm{q},\bm{g},\bm{h}\in\mathbb{X}^{n}$ and a matrix $\bm{B}\in\mathbb{X}^{n\times n}$, find regular vectors $\bm{x}\in\mathbb{X}^{n}$ that solve the problem
\begin{equation}
\begin{aligned}
&
\text{minimize}
&&
\bm{x}^{-}\bm{p}\oplus\bm{q}^{-}\bm{x},
\\
&
\text{subject to}
&&
\bm{B}\bm{x}\oplus\bm{g}
\leq
\bm{x},
\\
&
&&
\bm{x}
\leq
\bm{h}.
\end{aligned}
\label{P-xpqxBxxgxh}
\end{equation}

Note that the conditions in the problem can be represented in the equivalent form
\begin{align*}
&\bm{B}\bm{x}
\leq
\bm{x},
\\
&
\bm{g}
\leq
\bm{x}
\leq
\bm{h}.
\end{align*}

The solution of the problem follows an approach that is based on the general solution of linear inequalities obtained in \cite{Krivulin2006Solution,Krivulin2009Methods} and refined in \cite{Krivulin2013Amultidimensional}. We introduce an additional  parameter, which indicates the minimum value of the objective function. The problem is then reduced to the solving of a linear inequality with a matrix that depends on the parameter. We exploit the existence condition for solutions of the inequality to evaluate the parameter, and then take the general solution to the inequality as the solution to the initial problem.

\begin{theorem}[\cite{Krivulin2014Complete}]
\label{T-xpqxBxxgxh}
Let $\bm{x}$ be the general regular solution to problem \eqref{P-xpqxBxxgxh}, where $\bm{B}$ is a matrix such that $\mathop\mathrm{Tr}(\bm{B})\leq\mathbb{1}$, $\bm{p}$ is a nonzero vector, $\bm{q}$ and $\bm{h}$ are regular vectors, and $\bm{g}$ is a vector such that $\bm{h}^{-}\bm{B}^{\ast}\bm{g}\leq\mathbb{1}$. Then, the minimum in \eqref{P-xpqxBxxgxh} is equal to
\begin{equation*}
\theta
=
(\bm{q}^{-}\bm{B}^{\ast}\bm{p})^{1/2}
\oplus
\bm{h}^{-}\bm{B}^{\ast}\bm{p}\oplus\bm{q}^{-}\bm{B}^{\ast}\bm{g},
\end{equation*}
and the general solution is given by
\begin{equation*}
\bm{x}
=
\bm{B}^{\ast}\bm{u},
\qquad
\bm{g}\oplus\theta^{-1}\bm{p}
\leq
\bm{u}
\leq
((\bm{h}^{-}\oplus\theta^{-1}\bm{q}^{-})\bm{B}^{\ast})^{-}.
\end{equation*}
\end{theorem}

As a consequence, we obtain a complete solution to the problem, which was partially solved in \cite{Krivulin2012Anew}. Consider a variant of problem \eqref{P-xpqxBxxgxh} without the boundary constraints,
\begin{equation}
\begin{aligned}
&
\text{minimize}
&&
\bm{x}^{-}\bm{p}\oplus\bm{q}^{-}\bm{x},
\\
&
\text{subject to}
&&
\bm{B}\bm{x}
\leq
\bm{x}.
\end{aligned}
\label{P-xpqxBxx}
\end{equation}

\begin{corollary}[\cite{Krivulin2014Complete}]
\label{C-xpqxBxx}
Let $\bm{x}$ be the general regular solution to problem \eqref{P-xpqxBxx}, where $\bm{B}$ is a matrix such that $\mathop\mathrm{Tr}(\bm{B})\leq\mathbb{1}$, $\bm{p}$ is a nonzero vector and $\bm{q}$ is a regular vector.

Then, the minimum in \eqref{P-xpqxBxx} is equal to
\begin{equation*}
\theta
=
(\bm{q}^{-}\bm{B}^{\ast}\bm{p})^{1/2},
\end{equation*}
and the general solution is given by
\begin{equation*}
\bm{x}
=
\bm{B}^{\ast}\bm{u},
\qquad
\theta^{-1}\bm{p}
\leq
\bm{u}
\leq
\theta(\bm{q}^{-}\bm{B}^{\ast})^{-}.
\end{equation*}
\end{corollary}

\subsection{Problems with Span Seminorm}

The problems, in which the objective function is defined through the span seminorm, appeared in the context of job scheduling \cite{Butkovic2009Onsome,Tam2010Optimizing}. Minimization of the span seminorm solves scheduling problems in the just-in-time manufacturing. Maximization problems arise when the optimal schedule aims to spread the initiation or completion time of the jobs over the maximum possible time interval. 

The solution to the problem without constraints uses the evaluation of lower or upper bounds on the objective function. To solve the problems with linear equation and inequality constraints, we first obtain a general solution to the equation or inequality, and then substitute it into the objective function to get an unconstrained problem with known solution.

\subsubsection{Minimization Problems}

We start with an unconstrained problem that presents a rather extended version of minimization problem \eqref{P-minmax1AxAx1}. Given matrices $\bm{A},\bm{B}\in\mathbb{X}^{m\times n}$ and vectors $\bm{p},\bm{q}\in\mathbb{X}^{m}$, the problem is to find regular vectors $\bm{x}\in\mathbb{X}^{n}$ that
\begin{equation}
\begin{aligned}
&
\text{minimize}
&&
\bm{q}^{-}\bm{B}\bm{x}(\bm{A}\bm{x})^{-}\bm{p}.
\end{aligned}
\label{P-minqBxAxp}
\end{equation}

\begin{theorem}[\cite{Krivulin2013Explicit}]
\label{T-minqBxAxp}
Let $\bm{x}$ be a regular solution to problem \eqref{P-minqBxAxp}, where $\bm{A}$ is a row-regular matrix, $\bm{B}$ is a column-regular matrix, $\bm{p}$ is nonzero vector and $\bm{q}$ is regular vector.

Then, the minimum in \eqref{P-minqBxAxp} is equal to
\begin{equation*}
\Delta
=
(\bm{A}(\bm{q}^{-}\bm{B})^{-})^{-}\bm{p},
\end{equation*}
and attained at
\begin{equation*}
\bm{x}
=
\alpha(\bm{q}^{-}\bm{B})^{-},
\qquad
\alpha>\mathbb{0}.
\end{equation*}
\end{theorem}

We now examine particular cases of problem \eqref{P-minqBxAxp}. Suppose that $\bm{B}=\bm{A}=\bm{I}$ and $\bm{p}=\bm{q}=\mathbb{1}$, where $\mathbb{1}$ denotes the vector having all elements equal to $\mathbb{1}$, and consider the problem
\begin{equation*}
\begin{aligned}
&
\text{minimize}
&&
\mathbb{1}^{T}\bm{x}\bm{x}^{-}\mathbb{1}.
\end{aligned}
\end{equation*}

Application of Theorem~\ref{T-minqBxAxp} shows that the problem has the minimum $\Delta=\mathbb{1}$, which is attained at any vector $\bm{x}=\alpha\,\mathbb{1}$, $\alpha>\mathbb{0}$. 

If we assume that $\bm{A}=\bm{B}$ and $\bm{p}=\bm{q}=\mathbb{1}$, then we arrive at problem \eqref{P-minmax1AxAx1}. Theorem~\ref{T-minqBxAxp} allows us to obtain the solution to the last problem in a new compact vector form.

\begin{corollary}
\label{C-min1AxAx1}
Let $\bm{x}$ be a regular solution to problem \eqref{P-minmax1AxAx1}, where $\bm{A}$ is a regular matrix.

Then, the minimum in \eqref{P-minmax1AxAx1} is equal to
\begin{equation*}
\Delta
=
(\bm{A}(\mathbb{1}^{T}\bm{A})^{-})^{-}\mathbb{1},
\end{equation*}
and attained at
\begin{equation*}
\bm{x}
=
\alpha(\mathbb{1}^{T}\bm{A})^{-},
\qquad
\alpha>\mathbb{0}.
\end{equation*}
\end{corollary}

Consider a constrained problem: given matrices $\bm{B},\bm{C}\in\mathbb{X}^{n\times n}$, find regular vectors $\bm{x}\in\mathbb{X}^{n}$ such that
\begin{equation}
\begin{aligned}
&
\text{minimize}
&&
\mathbb{1}^{T}\bm{y}\bm{y}^{-}\mathbb{1},
\\
&
\text{subject to}
&&
\bm{C}\bm{x}
=
\bm{y},
\\
&
&&
\bm{D}\bm{x}
\leq
\bm{x}.
\end{aligned}
\label{P-1yy1CxyDxxI}
\end{equation}

\begin{theorem}[\cite{Krivulin2013Explicit}]
\label{T-1yy1CxyDxxI}
Let $\bm{x}$ be a regular solution to problem \eqref{P-1yy1CxyDxxI}, where $\bm{C}$ is a regular matrix and $\bm{D}$ is a matrix such that $\mathop\mathrm{Tr}(\bm{D})\leq\mathbb{1}$. Then, the minimum in \eqref{P-1yy1CxyDxxI} is equal to
$$
\Delta
=
(\bm{C}\bm{D}^{\ast}(\mathbb{1}^{T}\bm{C}\bm{D}^{\ast})^{-})^{-}\mathbb{1},
$$
and attained at
$$
\bm{x}
=
\alpha\bm{D}^{\ast}(\mathbb{1}^{T}\bm{C}\bm{D}^{\ast})^{-},
\qquad
\alpha
>
\mathbb{0}.
$$
\end{theorem}

\subsubsection{Maximization Problems}

We begin with an unconstrained problem. Given matrices $\bm{A}\in\mathbb{X}^{m\times n}$, $\bm{B}\in\mathbb{X}^{l\times n}$, and vectors $\bm{p}\in\mathbb{X}^{m}$, $\bm{q}\in\mathbb{X}^{l}$, we need to find regular vectors $\bm{x}\in\mathbb{X}^{n}$ that solve the problem
\begin{equation}
\begin{aligned}
&
\text{maximize}
&&
\bm{q}^{-}\bm{B}\bm{x}(\bm{A}\bm{x})^{-}\bm{p}.
\end{aligned}
\label{P-maxqBxAxp}
\end{equation}

The following result offers a complete solution to \eqref{P-maxqBxAxp} under fairly general conditions.
\begin{theorem}[\cite{Krivulin2013Amaximization}]
\label{T-maxqBxAxp}
Let $\bm{x}$ be the general regular solution to problem \eqref{P-maxqBxAxp}, where $\bm{A}$ is a matrix with regular columns, $\bm{B}$ is a column-regular matrix, $\bm{p}$ and $\bm{q}$ are regular vectors. Then, the minimum in \eqref{P-maxqBxAxp} is equal to
\begin{equation*}
\Delta
=
\bm{q}^{-}\bm{B}\bm{A}^{-}\bm{p},
\label{E-DeltaqBAp}
\end{equation*}
and the general solution $\bm{x}=(x_{i})$ has components that are given by
\begin{equation*}
\begin{aligned}
x_{k}
&=
\alpha\bm{a}_{k}^{-}\bm{p},
\\
x_{j}
&\leq
\alpha a_{sj}^{-1}p_{s},
\quad
j\ne k,
\end{aligned}
\label{I-xkalphaakp}
\end{equation*}
where $\alpha>\mathbb{0}$, and the indices $k$ and $s$ are defined by the conditions
$$
k
=
\arg\max_{1\leq i\leq n}\bm{q}^{-}\bm{b}_{i}\bm{a}_{i}^{-}\bm{p},
\qquad
s
=
\arg\max_{1\leq i\leq m}a_{ik}^{-1}p_{i}.
$$
\end{theorem}

Consider a special case of the problem. Assume that $\bm{p}=\bm{q}=\mathbb{1}$. Then, we have
$$
\mathbb{1}^{T}\bm{B}\bm{x}(\bm{A}\bm{x})^{-}\mathbb{1}
=
\|\bm{B}\bm{x}\|\|(\bm{A}\bm{x})^{-}\|,
\qquad
\mathbb{1}^{T}\bm{B}\bm{A}^{-}\mathbb{1}
=
\|\bm{B}\bm{A}^{-}\|.
$$

Under this assumption, problem \eqref{P-maxqBxAxp} becomes
\begin{equation}
\begin{aligned}
&
\text{maximize}
&&
\|\bm{B}\bm{x}\|\|(\bm{A}\bm{x})^{-}\|.
\end{aligned}
\label{P-max1BxAx1}
\end{equation}

Application of Theorem~\ref{T-maxqBxAxp} leads to the following result.
\begin{corollary}[\cite{Krivulin2013Amaximization}]
\label{C-max1BxAx1}
Let $\bm{x}$ be the general regular solution to problem \eqref{P-max1BxAx1}, where $\bm{A}$ is a matrix with regular columns and $\bm{B}$ is a column-regular matrix.

Then, the minimum in \eqref{P-max1BxAx1} is equal to 
\begin{equation*}
\Delta
=
\|\bm{B}\bm{A}^{-}\|,
\label{E-DeltaBA}
\end{equation*}
and the general solution $\bm{x}=(x_{i})$ has components that are given by
\begin{equation*}
\begin{aligned}
x_{k}
&=
\alpha\|\bm{a}_{k}^{-}\|,
\\
x_{j}
&\leq
\alpha a_{sj}^{-1},
\quad
j\ne k,
\end{aligned}
\label{I-xkalphabk}
\end{equation*}
where $\alpha>\mathbb{0}$, and the indices $k$ and $s$ are defined by the conditions
$$
k
=
\arg\max_{1\leq i\leq n}\|\bm{b}_{i}\|\|\bm{a}_{i}^{-}\|,
\qquad
s
=
\arg\max_{1\leq i\leq m}a_{ik}^{-1}.
$$
\end{corollary}

We now discuss the solution to constrained problems. Suppose that, given a matrix $\bm{C}\in\mathbb{X}^{n\times n}$, we need to solve the problem 
\begin{equation}
\begin{aligned}
&
\text{maximize}
&&
\bm{q}^{-}\bm{B}\bm{x}(\bm{A}\bm{x})^{-}\bm{p},
\\
&
\text{subject to}
&&
\bm{C}\bm{x}
\leq
\bm{x}.
\end{aligned}
\label{P-maxqBxAxpCxlex}
\end{equation}

It is known (see, e.g., \cite{Krivulin2009Methods}) that the inequality constraints have regular solutions if and only if $\mathop\mathrm{Tr}(\bm{C})\leq\mathbb{1}$. Under this condition, the general solution takes the form $\bm{x}=\bm{C}^{\ast}\bm{u}$, where $\bm{u}$ is any regular vector. Problem \eqref{P-maxqBxAxpCxlex} then reduces to the unconstrained problem
\begin{equation*}
\begin{aligned}
&
\text{maximize}
&&
\bm{q}^{-}\bm{B}\bm{C}^{\ast}\bm{u}(\bm{A}\bm{C}^{\ast}\bm{u})^{-}\bm{p}.
\end{aligned}
\label{P-maxqBCuACup-ast}
\end{equation*}

The last problem takes the form of \eqref{P-maxqBxAxp}, and thus is completely solved by Theorem~\ref{T-maxqBxAxp}. 

In a similar way, we obtain the solution to the problem with equality constraints
\begin{equation}
\begin{aligned}
&
\text{maximize}
&&
\bm{q}^{-}\bm{B}\bm{x}(\bm{A}\bm{x})^{-}\bm{p},
\\
&
\text{subject to}
&&
\bm{C}\bm{x}
=
\bm{x}.
\end{aligned}
\label{P-maxqBxAxpCxeq}
\end{equation}

Indeed, if the condition $\mathop\mathrm{Tr}(\bm{C})=\mathbb{1}$ holds, then the equation $\bm{C}\bm{x}=\bm{x}$ has the general solution $\bm{x}=\bm{C}^{+}\bm{u}$, where $\bm{u}$ is any regular vector \cite{Krivulin2009Methods}. Substitution of the solution into the objective function again leads to a problem in the form of \eqref{P-maxqBxAxp},
\begin{equation*}
\begin{aligned}
&
\text{maximize}
&&
\bm{q}^{-}\bm{B}\bm{C}^{+}\bm{u}(\bm{A}\bm{C}^{+}\bm{u})^{-}\bm{p}.
\end{aligned}
\label{P-maxqBCuACupplus}
\end{equation*}

\subsection{Problems with Evaluation of Spectral Radius}

First, we return to problem \eqref{P-minxAx}, which admits a solution based on the extremal property of the spectral radius. A complete solution to the problem was given in \cite{Krivulin2013Amultidimensional} as a consequence of solution to a more general optimization problem. A direct complete solution to problem \eqref{P-minxAx} is obtained in \cite{Krivulin2013Extremal} in the following form.
\begin{lemma}\label{L-minxAx}
Let $\bm{x}$ be the general regular solution to problem \eqref{P-minxAx}, where $\bm{A}$ is a matrix with spectral radius $\lambda>\mathbb{0}$.

Then, the minimum in \eqref{P-minxAx} is equal to $\lambda$, and the general solution is given by
$$
\bm{x}
=
(\lambda^{-1}\bm{A})^{\ast}\bm{u},
\qquad
\bm{u}\in\mathbb{X}^{n}.
$$
\end{lemma}

The proof of the statement is based on the derivation of a sharp lower bound for the objective function, accompanied by the development of the equation, which specifies the solution set for the problem. We then reduce the equation to an inequality and take the general solution to the inequality as a complete solution to the optimization problem.

Extensions of problem \eqref{P-minxAx}, which differ from \eqref{P-minxAx} by more general form of the objective function and constraints, were examined in \cite{Krivulin2012Acomplete,Krivulin2014Aconstrained}, where general solutions were obtained under various assumptions. 

Below, we offer complete direct solutions to certain new problems, which present further extensions of problem \eqref{P-minxAx}.

\subsubsection{An Unconstrained Problem}

Suppose that a matrix $\bm{A}\in\mathbb{X}^{n\times n}$, vectors $\bm{p},\bm{q}\in\mathbb{X}^{n}$, and a scalar $c\in\mathbb{X}$ are given. Consider the problem of finding regular vectors $\bm{x}\in\mathbb{X}^{n}$ that
\begin{equation}
\begin{aligned}
&
\text{minimize}
&&
\bm{x}^{-}\bm{A}\bm{x}\oplus\bm{x}^{-}\bm{p}\oplus\bm{q}^{-}\bm{x}\oplus c.
\end{aligned}
\label{P-xAxxpqxc}
\end{equation}

The problem is completely solved as follows.
\begin{theorem}[\cite{Krivulin2013Extremal}]
\label{T-xAxxpqxc}
Let $\bm{x}$ be the general regular solution to problem \eqref{P-xAxxpqxc}, where $\bm{A}$ is a matrix with spectral radius $\lambda>\mathbb{0}$ and $\bm{q}$ is a regular vector. Denote
\begin{equation*}
\mu
=
\lambda
\oplus
\bigoplus_{m=1}^{n}
(\bm{q}^{-}\bm{A}^{m-1}\bm{p})^{1/(m+1)}
\oplus
c.
\end{equation*}

Then, the minimum in \eqref{P-xAxxpqxc} is equal to $\mu$ and the general solution is given by
\begin{equation*}
\bm{x}
=
(\mu^{-1}\bm{A})^{\ast}\bm{u},
\qquad
\mu^{-1}\bm{p}
\leq
\bm{u}
\leq
\mu(\bm{q}^{-}(\mu^{-1}\bm{A})^{\ast})^{-}.
\end{equation*}
\end{theorem}

\subsubsection{Problems with Constraints}

Given matrices $\bm{A},\bm{B}\in\mathbb{X}^{n\times n}$, $\bm{C}\in\mathbb{X}^{m\times n}$, and vectors $\bm{g}\in\mathbb{X}^{n}$, $\bm{h}\in\mathbb{X}^{m}$, consider the problem to find regular vectors $\bm{x}\in\mathbb{X}^{n}$ that
\begin{equation}
\begin{aligned}
&
\text{minimize}
&&
\bm{x}^{-}\bm{A}\bm{x},
\\
&
\text{subject to}
&&
\bm{B}\bm{x}\oplus\bm{g}
\leq
\bm{x},
\\
&&&
\bm{C}
\bm{x}
\leq
\bm{h}.
\end{aligned}
\label{P-xAxBxgxCxh}
\end{equation}

\begin{theorem}[\cite{Krivulin2014Aconstrained}]
\label{T-xAxBxgxCxh}
Let $\bm{x}$ be the general regular solution to problem \eqref{P-xAxBxgxCxh}, where $\bm{A}$ is a matrix with spectral radius $\lambda>\mathbb{0}$, $\bm{B}$ is a matrix such that $\mathop\mathrm{Tr}(\bm{B})\leq\mathbb{1}$. Suppose that $\bm{C}$ is a column-regular matrix and $\bm{h}$ is a regular vector such that $\bm{h}^{-}\bm{C}\bm{B}^{\ast}\bm{g}\leq\mathbb{1}$.

Then, the minimum in \eqref{P-xAxBxgxCxh} is equal to
\begin{equation*}
\theta
=
\bigoplus_{k=1}^{n}\mathop{\bigoplus\hspace{0.0em}}_{0\leq i_{0}+
i_{1}+\cdots+i_{k}\leq n-k}\mathop\mathrm{tr}\nolimits^{1/k}(\bm{B}^{i_{0}}
(\bm{A}\bm{B}^{i_{1}}\cdots\bm{A}\bm{B}^{i_{k}})(\bm{I}\oplus\bm{g}\bm{h}^{-}\bm{C})),
\end{equation*}
and the general solution is given by
\begin{equation*}
\bm{x}
=
(\theta^{-1}\bm{A}\oplus\bm{B})^{\ast}\bm{u},
\qquad
\bm{g}
\leq
\bm{u}
\leq
(\bm{h}^{-}\bm{C}(\theta^{-1}\bm{A}\oplus\bm{B})^{\ast})^{-}.
\end{equation*}
\end{theorem}

In the case when $\bm{C}=\mathbb{0}$, problem \eqref{P-xAxBxgxCxh} takes the form 
\begin{equation}
\begin{aligned}
&
\text{minimize}
&&
\bm{x}^{-}\bm{A}\bm{x},
\\
&
\text{subject to}
&&
\bm{B}\bm{x}\oplus\bm{g}
\leq
\bm{x}.
\end{aligned}
\label{P-xAxBxgx}
\end{equation}

\begin{corollary}[\cite{Krivulin2013Amultidimensional,Krivulin2014Aconstrained}]
\label{C-xAxBxgx}
Let $\bm{x}$ be the general regular solution to problem \eqref{P-xAxBxgx}, where $\bm{A}$ is a matrix with spectral radius $\lambda>\mathbb{0}$. Then, the minimum in \eqref{P-xAxBxgx} is equal to
$$
\theta
=
\lambda
\oplus
\bigoplus_{k=1}^{n-1}\mathop{\bigoplus\hspace{1.2em}}_{1\leq i_{1}+\cdots+i_{k}\leq n-k}
\mathop\mathrm{tr}\nolimits^{1/k}(\bm{A}\bm{B}^{i_{1}}\cdots\bm{A}\bm{B}^{i_{k}}),
$$
and the general solution is given by
$$
\bm{x}
=
(\theta^{-1}\bm{A}\oplus\bm{B})^{\ast}\bm{u},
\qquad
\bm{u}
\geq
\bm{g}.
$$
\end{corollary}

Suppose that $\bm{B}=\mathbb{0}$ and $\bm{C}=\bm{I}$. In this case, problem \eqref{P-xAxBxgxCxh} becomes
\begin{equation}
\begin{aligned}
&
\text{minimize}
&&
\bm{x}^{-}\bm{A}\bm{x},
\\
&
\text{subject to}
&&
\bm{g}
\leq
\bm{x}
\leq
\bm{h}.
\end{aligned}
\label{P-xAxgxh}
\end{equation}

\begin{corollary}[\cite{Krivulin2014Aconstrained}]
\label{C-xAxgxh}
Let $\bm{x}$ be the general regular solution to problem \eqref{P-xAxgxh}, where $\bm{A}$ is a matrix with spectral radius $\lambda>\mathbb{0}$, $\bm{g}$ is a vector and $\bm{h}$ is a regular vector such that $\bm{h}^{-}\bm{g}\leq\mathbb{1}$. Then, the minimum in \eqref{P-xAxgxh} is equal to
$$
\theta
=
\lambda
\oplus
\bigoplus_{k=1}^{n}(\bm{h}^{-}\bm{A}^{k}\bm{g})^{1/k},
$$
and the general solution is given by
$$
\bm{x}
=
(\theta^{-1}\bm{A})^{\ast}\bm{u},
\qquad
\bm{g}
\leq
\bm{u}
\leq
(\bm{h}^{-}(\theta^{-1}\bm{A})^{\ast})^{-}.
$$
\end{corollary}

The next problem combines a special case of the objective function of problem \eqref{P-xAxxpqxc} with the constraints of problem \eqref{P-xAxBxgx}. Given matrices $\bm{A},\bm{B}\in\mathbb{X}^{n\times n}$ and vectors $\bm{p},\bm{g}\in\mathbb{X}^{n}$, we find regular vectors $\bm{x}\in\mathbb{X}^{n}$ that provide solutions to the problem 
\begin{equation}
\begin{aligned}
&
\text{minimize}
&&
\bm{x}^{-}\bm{A}\bm{x}\oplus\bm{x}^{-}\bm{p},
\\
&
\text{subject to}
&&
\bm{B}\bm{x}\oplus\bm{g}
\leq
\bm{x}.
\end{aligned}
\label{P-xAxxpBxgx}
\end{equation}

\begin{theorem}[\cite{Krivulin2013Extremal}]
\label{T-xAxxpBxgx}
Let $\bm{x}$ be the general regular solution to problem \eqref{P-xAxxpBxgx}, where $\bm{A}$ is a matrix with spectral radius $\lambda>\mathbb{0}$ and $\bm{B}$ is a matrix such that $\mathop\mathrm{Tr}(\bm{B})\leq\mathbb{1}$.

Then, the minimum in \eqref{P-xAxxpBxgx} is equal to
\begin{equation*}
\theta
=
\lambda
\oplus
\bigoplus_{k=1}^{n-1}\mathop{\bigoplus\hspace{1.2em}}_{1\leq i_{1}+
\cdots+i_{k}\leq n-k}\mathop\mathrm{tr}\nolimits^{1/k}(\bm{A}\bm{B}^{i_{1}}\cdots\bm{A}\bm{B}^{i_{k}}),
\end{equation*}
and the general solution is given by
$$
\bm{x}
=
(\theta^{-1}\bm{A}\oplus\bm{B})^{\ast}\bm{u},
\qquad
\bm{u}
\geq
\theta^{-1}\bm{p}\oplus\bm{g}.
$$
\end{theorem}

\subsection*{Acknowledgments}
The author thanks an editor for helpful comments and suggestions on earlier drafts of the manuscript.

\bibliographystyle{utphys}

\bibliography{Tropical_optimization_problems}

\end{document}